\documentclass[review]{elsarticle}
\usepackage{tikz}
\usepackage{lineno,hyperref}
\usepackage{mathtools}
\DeclarePairedDelimiter\abs{\lvert}{\rvert}
\pdfoutput=1
\modulolinenumbers[5]
\usepackage{float}
\RequirePackage{amsmath}[2016/11/05]
\usepackage[ruled,vlined]{algorithm2e}
\usepackage{natbib}
\journal{Engineering Applications of Artificial Intelligence}
\date{}
\usepackage{geometry}
\newgeometry{vmargin={15mm}, hmargin={12mm,17mm}}
\usepackage{caption}
\usepackage{subcaption}

\begin{document}

\begin{frontmatter}

\title{Sensitivity analysis using Physics-informed neural networks}

\author[1,2]{John M. Hanna}
\author[1]{Jos\'e V. Aguado}
\author[1]{Sebastien~Comas-Cardona}
\author[2]{Ramzi Askri}
\author[1]{Domenico Borzacchiello}

\address[1]{Nantes Université, Ecole Centrale Nantes, CNRS, GeM, UMR 6183, 1 Rue de la Noë, 44300 Nantes, France}
\address[2]{Nantes Université, IRT Jules Verne, 44340 Bouguenais, France}

\begin{abstract}

The goal of this paper is to provide a simple approach to perform local sensitivity analysis using Physics-informed neural networks (PINN). The main idea lies in adding a new term in the loss function that regularizes the solution in a small neighborhood near the nominal value of the parameter of interest. The added term represents the derivative of the loss function with respect to the parameter of interest. The result of this modification is a solution to the problem along with the derivative of the solution with respect to the parameter of interest (the sensitivity). We call the new technique SA-PNN which stands for sensitivity analysis in PINN. The effectiveness of the technique is shown using four examples: the first one is a simple one-dimensional advection-diffusion problem to show the methodology, the second is a two-dimensional Poisson's problem with nine parameters of interest, and the third and fourth examples are one and two-dimensional transient two-phase flow in porous media problem.
\end{abstract}


\end{frontmatter}


\section{Introduction}

\subsection{Motivation and outline}

Engineering applications often involve modeling complex physical phenomena that are governed by partial differential equations (PDEs). Approximate solutions to the PDEs are obtained using numerical techniques such as finite elements or finite volumes. The solutions rely on input parameters that characterize material or chemical properties, as well as initial and boundary conditions. Understanding the impact of small variations or errors in these parameters on the PDE solution is a challenging task due to the inherent complexity of the equations involved. This exploration of how the solution of PDEs responds to perturbations in input parameters is known as sensitivity analysis, which is represented quantitatively as the derivative of the solution with respect to the input parameters. Applications of sensitivity analysis include aerodynamic optimization \cite{sensitivity_aerodynamic_rev}, shape optimization in solid mechanics \cite{sensitivity_shape_opti_solid}, injection molding \cite{sensitivity_injection_molding} and biomedical applications \cite{sensitivity_tumor, sensitivity_brain}.

Conceptually, the simplest method to perform sensitivity analysis is the finite difference method. It involves iteratively solving the system using any numerical method while altering one input parameter at a time while keeping the others constant. The sensitivity can be estimated by employing finite differences with the solutions obtained from different values of the modified parameter. The technique is easy to implement, however, one issue lies in the choice of the step size \cite{sensitivity_fd_step_size}. It has to be chosen such that the truncation error is minimized so that accurate derivative estimation is obtained, which can be done by reducing the step size and also minimizing the computer round-off error, which is increased by reducing the step size. To find a good trade-off value, a trial and error method is usually done which requires several system evaluations. Due to that reason, the computational effort grows for each calculation of the solution derivative with respect to one input parameter of interest \cite{sensitivity_diffusion_reaction}.

The discrete direct method is another technique to perform sensitivity analysis which can be more beneficial than the finite difference technique in some cases \cite{book_cahpter}. It is based on solving a differentiated version of the original discretized system of equations with respect to the input parameter of interest. It requires having a discrete solution to the system, the derivative of the discretized stiffness matrix, and the force vector with respect to the input parameter of interest. In simple problems, analytical expressions can exist and the sensitivity can be efficiently obtained. However, in most cases, analytical expressions to these derivatives do not exist and approximations are obtained using finite differences. It was shown by \cite{discrete_problem} that using an approximation to the stiffness matrix derivative usually causes serious accuracy issues.  

The continuum direct method has a similar general methodology as the discrete direct method but it deals with the continuous form of the equations instead \cite{book_cahpter}. First, the differentiated form of the continuous equations is derived. This continuum sensitivity equation version can then be solved using numerical techniques. The major advantage of this approach is that the derived sensitivity equation is independent of the chosen numerical technique. The continuum and discrete approaches can lead to the same sensitivity solution if the same numerical technique is chosen, however, in the case of shape design, the two methods can lead to different results \cite{continuum_direcetre}. The method has similar drawbacks as the discrete approach when it comes to approximating terms representing derivatives with respect to the parameter of interest.

The adjoint method for both the discrete and continuum approaches provides big advantages to using the direct methods \cite{adjoint_review}. It is based on solving an adjoint system of equations, typically derived from the original system, which involves solving a linearized system of equations. Afterwards, the sensitivity is obtained from an equation substitution not based on solving a linear system. The main advantage lies in the independency of the adjoint system of the input parameters of interest so the linear system of equations can only be solved once \cite{adjoint_review}. Then, sensitivity can be obtained with respect to any number of input parameters of interest with low additional cost. For that reason, the adjoint method has become popular to perform sensitivity analysis, particularly in Computational Fluid Dynamics applications \cite{sensitivity_CFD_1}. Despite its widespread success, certain challenges persist. Notably, issues arise regarding the differentiability of solutions in the presence of discontinuities, which commonly occur in scenarios involving shock waves or two-phase flow which can lead to instabilities \cite{adjoint_two_phase}.

Machine learning-based techniques have experienced significant growth in addressing problems governed by partial differential equations (PDEs). Among these techniques, Physics-informed neural networks (PINN) have emerged as a rapidly expanding field focused on solving forward and inverse problems associated with PDEs \cite{pinn_main}. The appeal of PINN lies in its ability to circumvent the need for large datasets, as it leverages the underlying physics described by the PDEs to regularize the model. PINN's strength lies in utilizing feed-forward neural networks, which serve as universal approximators \cite{universal_approx}, and benefiting from recent advancements in automatic differentiation capabilities \cite{automatic_diff} to facilitate derivative calculations. Consequently, PINN has found applications in diverse areas, including solid mechanics \cite{pinn_solid}, fluid mechanics \cite{pinn_fluid_review, pinn_hidden_fluid}, additive manufacturing \cite{pinn_additive}, two-phase flow in porous media \cite{mine}, and numerous other domains \cite{xpinn, conservative_pinn, fractional_pinn, ultrasound_pinn, composite_pinn, surrogate_pinn, subsurface_pinn, fracture_pinn}.

This article utilizes PINN as the foundational framework to develop a novel sensitivity analysis technique. PINN framework is chosen due to the good approximation capabilities of neural networks, which were also shown in previous studies to accurately predict discontinuous fields \cite{mine} Moreover, it is simple and fast to get required derivatives using automatic differentiation. The primary contributions of the article can be summarized as follows:

\begin{itemize}
    \item Introducing a new technique, termed SA-PINN standing for sensitivity analysis in PINN, that leverages the PINN framework to perform sensitivity analysis. This is done by adding an extra term in the loss function that regularizes the solution within a neighborhood of nominal values of parameters of interest. 
    \item Demonstrating the simplicity and effectiveness of obtaining sensitivities simultaneously for multiple input parameters.
    \item Highlighting the capability of SA-PINN to accurately capture sensitivities in problems featuring discontinuities, including scenarios with moving boundaries or sharp gradients.
\end{itemize}

The structure of the paper is as follows. In section \ref{subsection:related}, the related work in the literature is discussed. In Chapter \ref{section:SA-pinn}, a brief introduction to PINN is provided, followed by an explanation of the SA-PINN technique. Chapter \ref{section:models} presents four problems investigated in the paper: a one-dimensional (1D) advection-diffusion problem, a two-dimensional (2D) Poisson's problem with multiple parameters of interest, and a 1D and 2D unsteady two-phase flow in porous media problem. Finally, Chapter \ref{section:conclusion} offers a discussion and conclusion to the work.

\subsection{Related work}
\label{subsection:related}

Parametric models can be utilized to perform local sensitivity analysis for parameters of interest. This is done by taking the derivative of the model output with respect to the input parameter to the network with automatic differentiation. Parametric models were built with PINN in several applications involving magnetic problems \cite{magnetic_pinn}, solid mechanics \cite{solid_pinn}, and fluid mechanics \cite{geometry_pinn, navier_parametric}. Building such models can be computationally expensive because collocation points need to be sampled in the parametric space. This is one of the main obstacles this paper is tackling.

A similar approach to the one followed in this paper was adopted to improve PINN solutions \cite{gpinn}. In that previous work, the authors suggested adding an extra term that represents the derivative of the residual with respect to spatial coordinates which can be easily obtained thanks to automatic differentiation. By this it can be ensured that the residual is minimized at collocation points and their neighborhood thus needs fewer points to have accurate solutions. Another work that utilizes a similar approach was developed for a physics-informed DeepONet framework \cite{deeponet}. In this paper, we follow a similar approach but in the parametric space by adding a term representing the derivative of the residual with respect to the parameter of interest.

Other works utilize the power of automatic differentiation and the differentiability of the neural network models to solve inverse problems. This has been done in several applications involving high-speed flows \cite{high_speed_pinn}, unsteady groundwater flows \cite{pinn_ground_inverse}, and two-phase flow in porous media \cite{hanna2024self}.

\section{Sensitivity Analysis-PINN (SA-PINN)}
\label{section:SA-pinn}

We consider a PDE of the form

\begin{equation}
    u_t + \mu L(u) = 0, \ \ \ \mathbf{x} \in \Omega,\ t \in [0, T]
\end{equation}

where $u_t$ is the time derivative, $L$ a general differential operator and $\mu$ a material parameter. Initial and boundary conditions for the problems are defined as

\begin{gather}
    u(0, \mathbf{x}) = u_0\\
    u(t, \mathbf{x_D}) = u_D\\
    B(u(t, \mathbf{x_N})) = f(\mathbf{x_N})
\end{gather}

where $B$ is a differential operator, $\mathbf{x_D}$ the boundary where Dirichlet boundary condition is enforced and $\mathbf{x_N}$ the boundary where Neumann boundary condition is applied.

 \subsection{PINN}

To address the problem using PINN, the first step is to select the approximation space, which in this case is a feed-forward neural network. Automatic differentiation is utilized to form a combined loss function, which includes the residual of the PDE at specific spatiotemporal collocation points as well as the error in enforcing the initial and boundary conditions. The solution is obtained by updating the weights and biases of the neural network to minimize the loss function using optimization algorithms such as stochastic gradient descent, Adam, or BFGS.

The loss function can be expressed as a weighted sum of several terms which reads as:

\begin{equation}
Loss_{PINN} = \lambda_0\ loss_0 + \lambda_D\ loss_D + \lambda_N\ loss_N + \lambda_{f}\ loss_{f}
\end{equation}

where $\lambda_i$ are the weights assigned to each loss term, playing a crucial role in the optimization process. The choice of these weights is still an open research topic in the PINN research community \cite{self_weights, self_weights2, self_weights3}. A common way to choose these weights is with a trial and error process, which can be tedious if the number of the loss terms is big. One way to simplify this process is to non-dimensionalize the equations to ensure that the different loss terms have comparable magnitudes \cite{weighting_1, weighting_2, weighting_3}. In our examples, we choose our test problems to have domains or parameters that correspond to the non-dimensionalized form of the equations in addition to some quick manual tuning if needed through trial and error. The individual loss terms are defined as follows:

\begin{gather}
    loss_0 = \frac{1}{N_0}\sum_{i=1}^{N_0} r_0^2(t_0^i, \mathbf{x}_0^i) =\frac{1}{N_0}\sum_{i=1}^{N_0} ||u(t_0^i, \mathbf{x}_0^i) - u_0^i||^2,\\[2ex]
    loss_D = \frac{1}{N_D}\sum_{i=1}^{N_D} r_D^2(t_D^i, \mathbf{x}_D^i) =\frac{1}{N_D}\sum_{i=1}^{N_D} ||u(t_D^i, \mathbf{x}_D^i) - u_D^i||^2,\\[2ex]
    loss_N = \frac{1}{N_N}\sum_{i=1}^{N_N} r_N^2(t_N^i, \mathbf{x}_N^i) = \frac{1}{N_N}\sum_{i=1}^{N_N} ||B(u(t_N^i, \mathbf{x}_N^i))-f_N^i||^2,\\[2ex]
    loss_{f} = \frac{1}{N_{f}}\sum_{i=1}^{N_f} ||r(t_f^i, \mathbf{x}_f^i)||^2= \frac{1}{N_f}\sum_{i=1}^{N_f} ||u_t + \mu L(u)||_{(t_f^i, \mathbf{x}_f^i)}^2.\\
\end{gather}

where $\{t_0^i, \mathbf{x_0^i}, u_0^i\}_{i=1}^{N_0}$ denote the training data for the initial condition, $\{t_D^i, \mathbf{x_D^i}, u_D^i\}_{i=1}^{N_D}$ specify the data for the Dirichlet boundary condition, $\{t_N^i, \mathbf{x_N^i}, f_N^i\}_{i=1}^{N_N}$ represent the data for the Neumann boundary condition, and $\{t_f^i, \mathbf{x_f^i}\}_{i=1}^{N_f}$ specify the collocation points where the residual of the PDE is being minimized. The subscripts $0$, $D$, $N$, and $f$ refer to the initial condition, Dirichlet condition, Neumann condition, and PDE residual, respectively.

 \subsection{SA-PINN}
 
The primary objective of PINN is to find a solution that minimizes the residual of the PDE within the spatiotemporal domain represented by the collocation points while satisfying the initial and boundary conditions. This yields a solution, denoted as $\hat{u}(t, \mathbf{x}; \hat{\mu})$, to the PDE at a specific value $\hat{\mu}$. When performing sensitivity analysis with $\mu$ as the input parameter of interest, we aim to not only determine the solution $\hat{u}$ but also calculate its derivative with respect to $\mu$, referred to as the sensitivity, at a given nominal value $\hat{\mu}$.

To obtain the sensitivity using PINN, a parametric or meta-model can be constructed. This involves modifying the structure of the neural network to incorporate an additional input, which is the parameter of interest $\mu$. Collocation points are then added in the spatiotemporal-parametric space. The residual is minimized across the entire parametric domain while respecting the initial and boundary conditions. Subsequently, the derivative of the solution with respect to $\mu$ can be directly obtained through automatic differentiation.

However, a challenge arises when building these parametric models: the number of collocation points grows exponentially as the number of parameters of interest increases. Consequently, computational intractability arises when dealing with multiple input parameters of interest.

To address this challenge, an alternative approach is proposed. Similarly to building parametric models, the structure of the network is modified to incorporate extra inputs representing the input parameters of interest. The collocation points are kept in the spatiotemporal domain and no points are added in the parametric space. Instead of solely minimizing the loss function, representing the residual of the PDE and the conditions, the derivative of the loss function with respect to the parameter of interest is also minimized. The modified loss function will then be formulated as the sum of the residual, the derivative of the residual with respect to $\mu$, and terms related to satisfying the initial and boundary conditions. By employing this technique, the solution can be accurately determined within a small neighborhood of $\hat{\mu}$, facilitating the calculation of sensitivity.

In summary, the SA-PINN technique can be outlined as follows:
 
\begin{itemize}
    \item Choose the neural network to have inputs related to space, time, and parameters of interest.
    \item Sample the collocation points only in space and time.
    \item Create the loss function having terms related to PDE residual, the residual derivative with respect to the parameter of interest, and the terms related to the initial and boundary conditions.
\end{itemize}

The modified loss function will then be

\begin{equation}
    Loss_{SA-PINN} = Loss_{PINN}  + Loss_{SA}
\end{equation}

where

\begin{equation}
    Loss_{SA} = \lambda_{0\mu}\ loss_{0\mu} + \lambda_{D\mu}\ loss_{D\mu} +\lambda_{N\mu}\ loss_{N\mu}+\lambda_{1\mu}\ loss_{r\mu}.
\end{equation}

and

\begin{equation}
    loss_{0\mu} = \frac{1}{N_0}\sum_{i=1}^{N_0} \abs*{\abs*{\frac{\partial r_0}{\partial \mu}}}^2_{(t_0^i, \mathbf{x}_0^i, \hat{\mu})},
\end{equation}

\begin{equation}
    loss_{D\mu} = \frac{1}{N_D}\sum_{i=1}^{N_D} \abs*{\abs*{\frac{\partial r_D}{\partial \mu}}}^2_{(t_D^i, \mathbf{x}_D^i, \hat{\mu})},
\end{equation}

\begin{equation}
    loss_{N\mu} = \frac{1}{N_N}\sum_{i=1}^{N_N} \abs*{\abs*{\frac{\partial r_N}{\partial \mu}}}^2_{(t_N^i, \mathbf{x}_N^i, \hat{\mu})},
\end{equation}

\begin{equation}
    loss_{f\mu} = \frac{1}{N_{f}}\sum_{i=1}^{N_f} \abs*{\abs*{\frac{\partial r}{\partial \mu}}}^2_{(t_f^i, \mathbf{x}_f^i, \hat{\mu})}.
\end{equation}

\begin{figure}[H]
    \centering
    \includegraphics[width=0.8\textwidth]{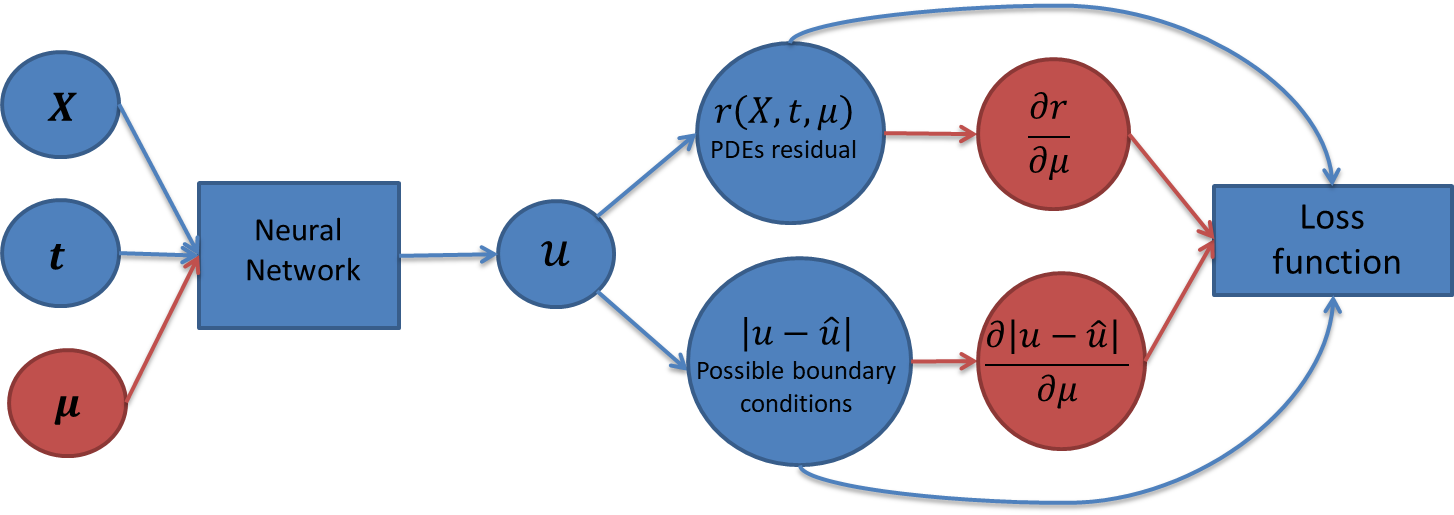}
  \caption{Diagram explaining the methodology of SA-PINN. The solution of the PDE is approximated by a neural network that takes an input of space $\mathbf{X}$, time $t$, and the parameter of interest $\mu$. Automatic differentiation is used to obtain the residual of the PDE. An extra step is performed to obtain the derivative of the residual with respect to the parameter of interest. The total loss function is composed of the PDE residual and boundary conditions term, as well as the derivative of the residual and boundary conditions term with respect to the parameter of interest.}
  \label{fig:sa_pinn}
\end{figure}

Figure~\ref{fig:sa_pinn} shows a diagram that summarizes the methodology of SA-PINN. The parts in orange are the added parts from classical PINN. The $u-\hat{u}$ term represents the mismatch of the solution from the initial and boundary conditions. It must be noted that we sample the collocation points only in space and time, but the points have another coordinate $\mu$ and all have a nominal value $\hat{\mu}$.

\section{Numerical examples}\label{section:models}

In this section, we introduce three different numerical examples that are used to show the effectiveness of the technique. Following the previous works \cite{mine, hanna2024self}, the chosen neural networks have 5 hidden layers with 20 neurons in each layer. The hyperbolic tangent activation function is applied to the hidden layers. For training, Adam optimizer is used first to escape local minima regions followed by the BFGS optimizer which is powerful enough to quickly reach the minimum in fewer iterations but requires a higher computational cost per iteration than Adam.

In all the examples, PINN refers to vanilla PINN where an extra input is added to represent the parameter of interest but with no points added to the parametric space and no modification to the loss function. This is chosen to show that vanilla PINN trained only on the nominal values does not generalize to other values of the parameter.

\subsection{1D diffusion-advection equation}

The first example is a steady one-dimensional diffusion-advection equation where we would like to study the effect of perturbations in the diffusion term $\epsilon$ on the solution. The strong form of the problem can be written as follows:
\begin{equation}
    \begin{gathered}
    \epsilon\ u_{xx} - u_x + 1 = 0,\ \ \ x \in [0,1],\\
    u(0) = 1,\ \ u(1) = 3\\
\end{gathered}
\end{equation}

The chosen nominal value for $\epsilon$ is $0.1$. $u_{xx}$ and $u_x$ are respectively the second and first-order derivatives of the solution $u$. The loss function is written as:

\begin{equation}
    Loss = \lambda_f\ loss_f + \lambda_b\ loss_b + \lambda_{f_\epsilon}\ loss_{r_\epsilon} + \lambda_{b_\epsilon}\ loss_{b_\epsilon}, 
\end{equation}

where 

\begin{equation}
    loss_{f} = \frac{1}{N_{f}}\sum_{i=1}^{N_f} ||r(x_f^i, \hat{\epsilon})||^2 =  \frac{1}{N_f}\sum_{i=1}^{N_f} ||\hat{\epsilon}\ u_{xx} - u_x + 1||_{x_f^i}^2,
\end{equation}
\begin{equation}
    loss_b = \frac{1}{N_b}\sum_{i=1}^{N_b} r_b^2(x_b^i, \hat{\epsilon}) =\frac{1}{N_b}\sum_{i=1}^{N_b} ||u(x_b^i, \hat{\epsilon}) - u_b^i||^2,
\end{equation}
\begin{equation}
    loss_{f_\epsilon} = \frac{1}{N_{f}}\sum_{i=1}^{N_f}\abs*{\abs*{\frac{\partial r}{\partial \epsilon}}}^2_{(x_f^i, \hat{\epsilon})},
\end{equation}
\begin{equation}
    loss_{b_\epsilon} = \frac{1}{N_b}\sum_{i=1}^{N_b}\abs*{\abs*{\frac{\partial r}{\partial \epsilon}}}^2_{(x_b^i, \hat{\epsilon})}.
\end{equation}

The weights for the different terms in the loss function are set to $1$ for the original PINN terms and $0.1$ for the added sensitivity terms. Adam optimizer is first used for 1000 iterations following, then the BFGS optimizer is used. 100 equally spaced collocation points are chosen for the training process.

\begin{figure}[H]
    \centering
    \includegraphics[width=0.7\textwidth]{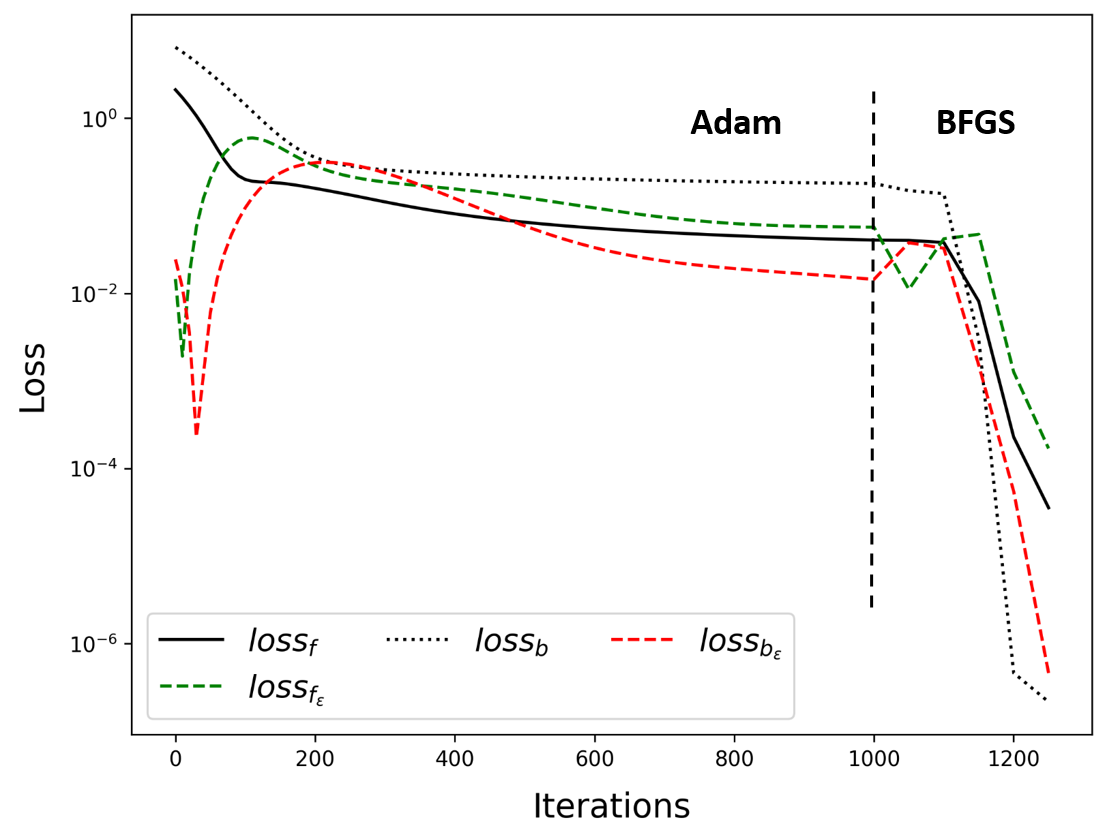}
  \caption{Loss function contribution for the 1D diffusion-advection example vs. the number of iterations.}\label{fig:loss_terms_1d}
\end{figure}

The different terms of the loss function are plotted against the iterations of the minimization algorithm in figure~\ref{fig:loss_terms_1d}.

\begin{figure}[H]
    \centering
    \includegraphics[width=0.7\textwidth]{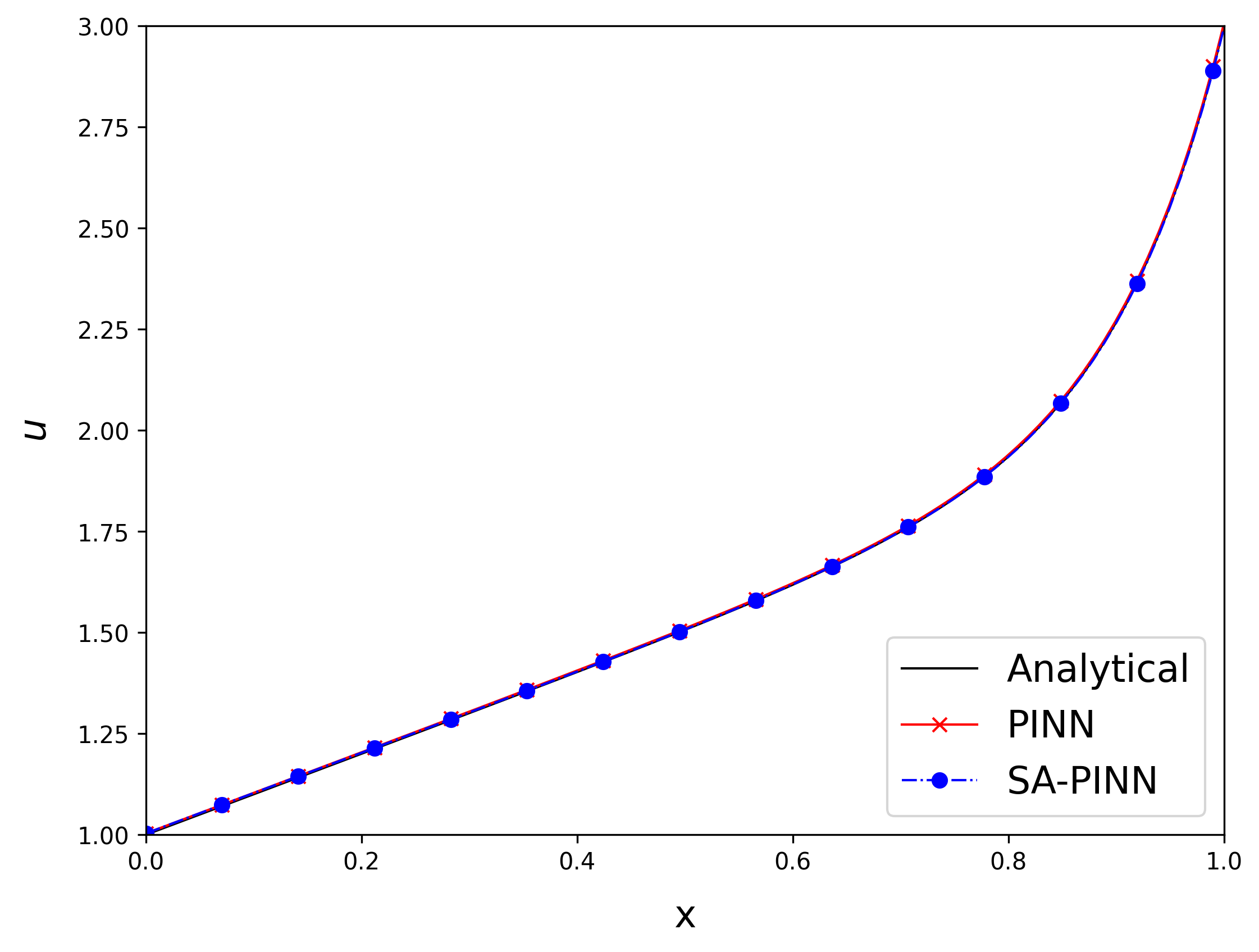}
  \caption{Solution $u$ at $\epsilon=0.1$ using PINN and SA-PINN along with the analytical solution of the 1D advection-diffusion problem.}\label{fig:1D_sol}
\end{figure}

The solution $u$ using PINN and SA-PINN is shown in figure~\ref{fig:1D_sol} along with the analytical solution for $\epsilon=0.1$.

\begin{figure}[H]
    \centering
    \includegraphics[width=0.7\textwidth]{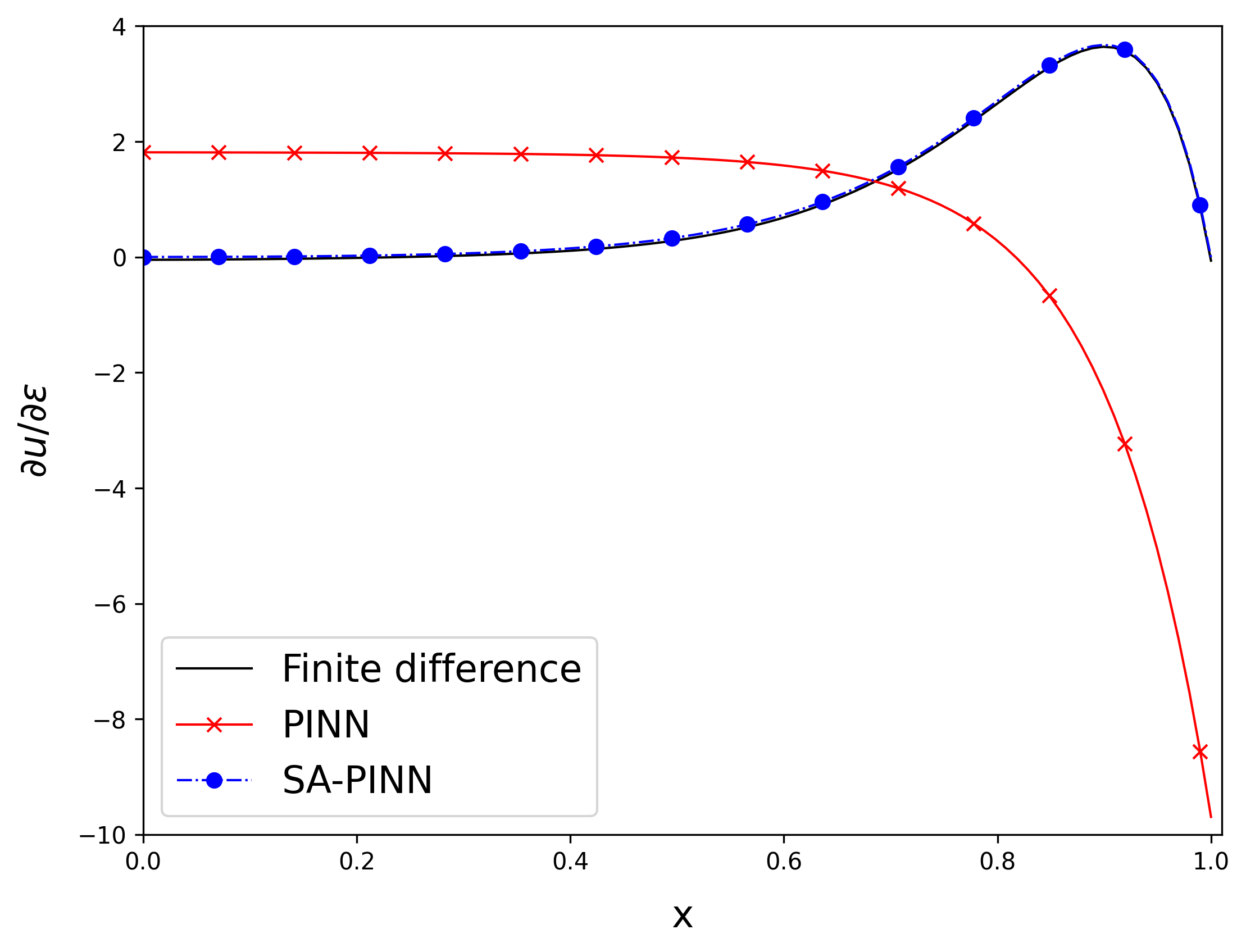}
  \caption{$\frac{\partial u}{\partial \epsilon}$ at $\epsilon=0.1$ using PINN and SA-PINN along with the finite difference solution of the 1D advection-diffusion problem.}
  \label{fig:1D_der}
\end{figure}

From figure~\ref{fig:1D_sol}, we can see that PINN and SA-PINN accurately capture the analytical solution to the problem. The derivative of the solution with respect to $\epsilon$ at $\epsilon=0.1$ is shown in figure~\ref{fig:1D_der}.

\begin{figure}[H]
    \centering
    \includegraphics[width=0.7\textwidth]{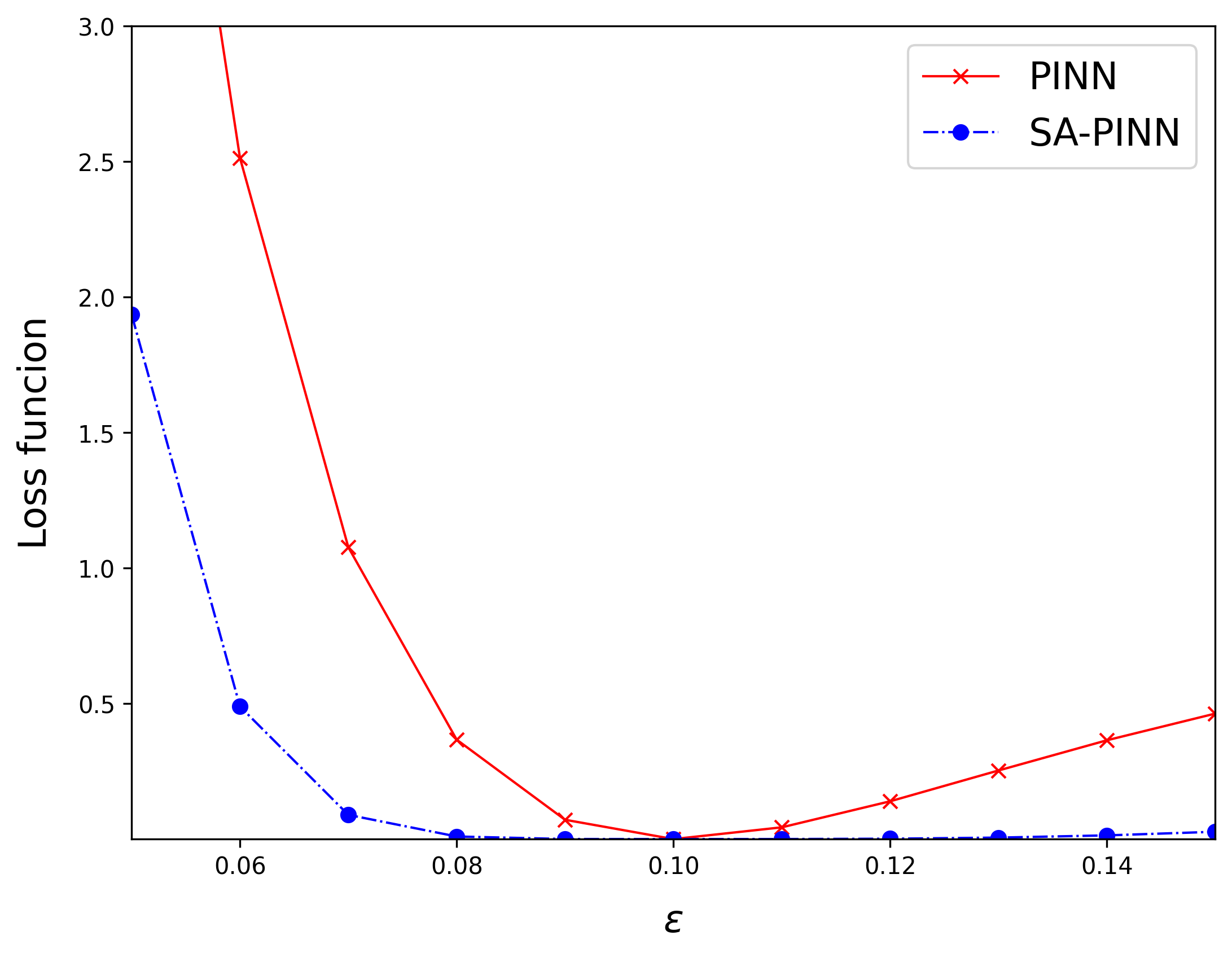}
  \caption{Loss term $loss_f$ for different $\epsilon$ values using PINN and SA-PINN of the 1D advection-diffusion problem.}
  \label{fig:1D_loss}
\end{figure}

The reference finite difference solution in figure~\ref{fig:1D_der} is obtained by obtaining different PINN solutions near $\epsilon=0.1$ and then calculating the derivative. We can see that classical PINN fails to predict the derivative, while, SA-PINN accurately predicts the derivative due to the added regularization term in the loss function. The loss term $loss_f$ is plotted after the training for different values of $\epsilon$ in figure~\ref{fig:1D_loss} for PINN and SA-PINN.

As seen in figure~\ref{fig:1D_loss}, SA-PINN has the effect of greatly flattening the loss curve in a neighborhood ($+/- 30\%$) near the nominal value of $\epsilon=0.1$. This leads to better solutions that PINN in the neighborhood and accurate derivative calculation at $\epsilon=0.1$.


\subsection{2D Poisson's problem}

\begin{figure}[H]
    \centering
\tikzset{every picture/.style={line width=0.75pt}} 

\begin{tikzpicture}[x=0.75pt,y=0.75pt,yscale=-1,xscale=1]

\draw  [fill={rgb, 255:red, 65; green, 117; blue, 5 }  ,fill opacity=1 ] (195.86,41) -- (358,41) -- (358,203.14) -- (195.86,203.14) -- cycle ;
\draw    (251,42.14) -- (250,203.14) ;
\draw    (305,42.14) -- (304,203.14) ;
\draw    (196,96.14) -- (358,97.14) ;
\draw    (196,148.14) -- (358,149.14) ;

\draw (213,57) node [anchor=north west][inner sep=0.75pt]  [font=\large]  {$k_{1}$};
\draw (323,57) node [anchor=north west][inner sep=0.75pt]  [font=\large]  {$k_{3}$};
\draw (267,57) node [anchor=north west][inner sep=0.75pt]  [font=\large]  {$k_{2}$};
\draw (215,111) node [anchor=north west][inner sep=0.75pt]  [font=\large]  {$k_{4}$};
\draw (320,160) node [anchor=north west][inner sep=0.75pt]  [font=\large]  {$k_{9}$};
\draw (320,111) node [anchor=north west][inner sep=0.75pt]  [font=\large]  {$k_{6}$};
\draw (268,112) node [anchor=north west][inner sep=0.75pt]  [font=\large]  {$k_{5}$};
\draw (270,162) node [anchor=north west][inner sep=0.75pt]  [font=\large]  {$k_{8}$};
\draw (214,162) node [anchor=north west][inner sep=0.75pt]  [font=\large]  {$k_{7}$};
\end{tikzpicture}

    \caption{2D Poisson's problem domain}
    \label{fig:poisson_domain}
\end{figure}
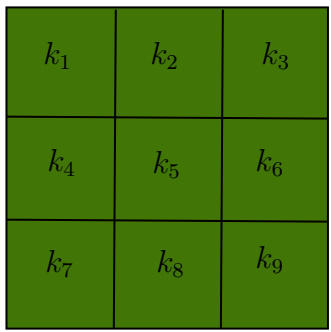

The next example is a 2-dimensional Poisson's problem where we have multiple parameters to study their effect on the solution. The domain is shown in figure~\ref{fig:poisson_domain} where there exist 9 subdomains each having different diffusivity values.

The strong form of the problem can be written as:

\begin{equation}
    \begin{gathered}
    k\  \Delta u = -1,\ \ \ in\ \Omega,\\
    u = 0, \ \ \ \ on\  \partial\Omega
\end{gathered}
\end{equation}

where $\Omega$ is a square with unit sides and $k$ is the diffusivity. The 9 subdomains have equal areas. The diffusivity parameters have the same nominal value which is $1$; $\hat{k_1}=\hat{k_2}=...=\hat{k_9}=1$.

The approximation space is chosen such that the boundary conditions are satisfied automatically. The approximation reads as follows:

\begin{equation}
    \hat{u} = x y  (x-1)  (y-1)\ u_{NN}(x, y, k_1,...,k_9)
\end{equation}

where $u_{NN}$ is a neural network with inputs $x$ and $y$ and the diffusivity parameters. The full loss function will then reads as:

\begin{equation}
    Loss = \lambda_f\ loss_f + \sum_{i=1}^{9}\lambda_i\ loss_{f_{k_i}},
    \label{eq:loss_2d_poisson}
    \end{equation}

where $\lambda_f=1$, all values of $\lambda_i$ are set to $0.1$,

\begin{equation}
    loss_{f} = \frac{1}{N_{f}}\sum_{i=1}^{N_f} ||r(x_f^i, \hat{k})||^2 =  \frac{1}{N_f}\sum_{i=1}^{N_f} ||\hat{k}\  \Delta u +1||_{x_f^i}^2,
\end{equation}

and

\begin{equation}
    loss_{f_{k_i}} = \frac{1}{N_{f}}\sum_{i=1}^{N_f}\abs*{\abs*{\frac{\partial r}{\partial k_i}}}^2_{(x_f^i, \hat{k_i})}.
\end{equation}

\begin{figure}[H]
    \centering
    \includegraphics[width=0.7\textwidth]{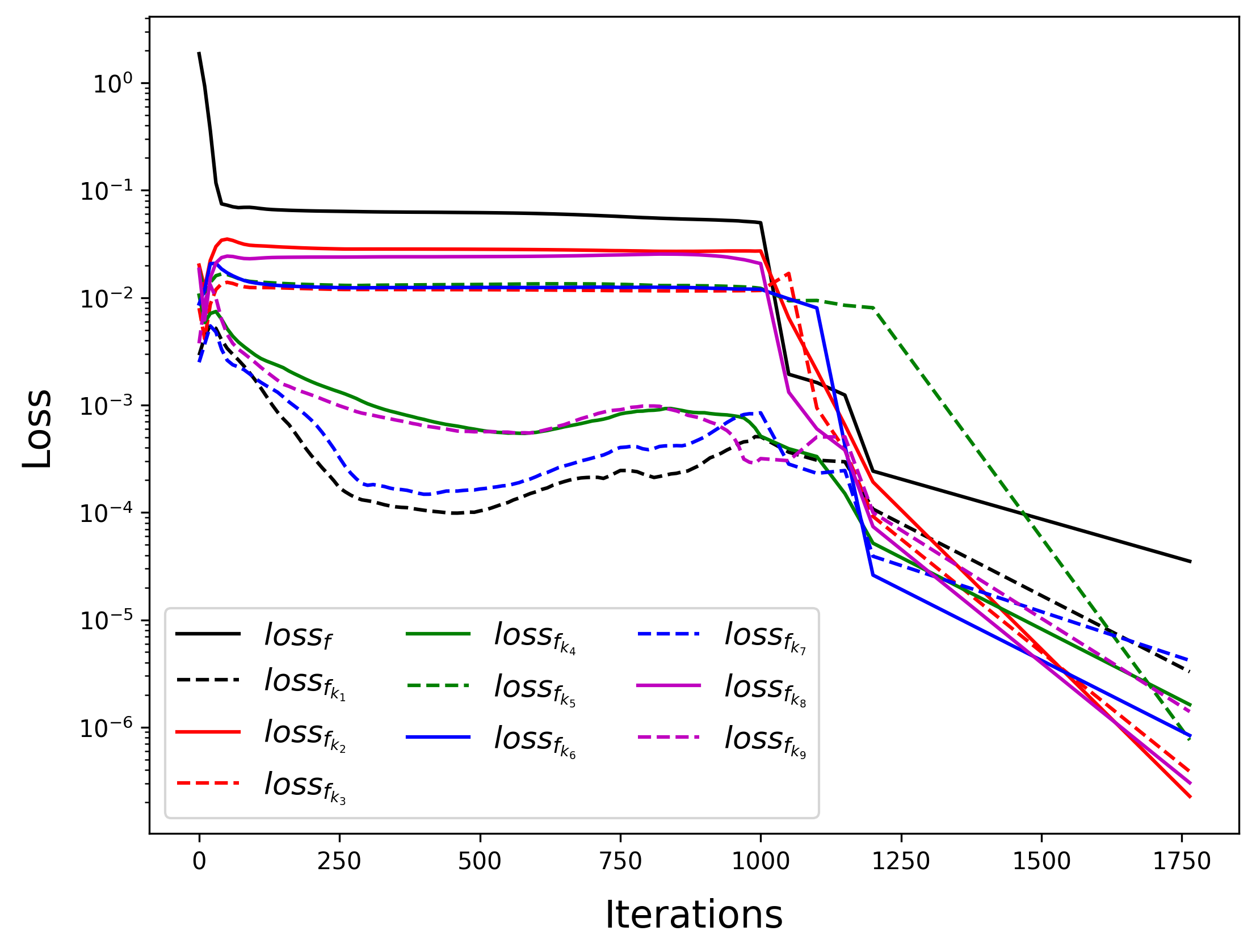}
  \caption{Loss function contribution for the 2D Poisson's example vs. the number of iterations.}\label{fig:loss_terms_2d}
\end{figure}

1000 collocation points sampled using the Latin hypercube sampling strategy are used for the training. The different loss terms are plotted against the iterations in figure~\ref{fig:loss_terms_2d}.

A finite difference code of the problem is developed and mesh sensitivity is studied to make sure that the solution is converged. The finite difference solution is used as the ground truth in this case

\begin{figure}[H]
     \centering
     \begin{subfigure}[b]{0.49\textwidth}
         \centering
    \includegraphics[width=1.1\textwidth]{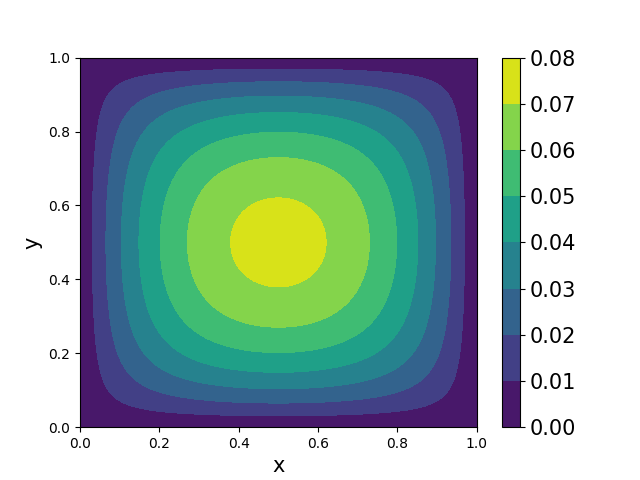}
     \end{subfigure}
     \hfill
     \begin{subfigure}[b]{0.49\textwidth}
         \centering
    \includegraphics[width=1.1\textwidth]{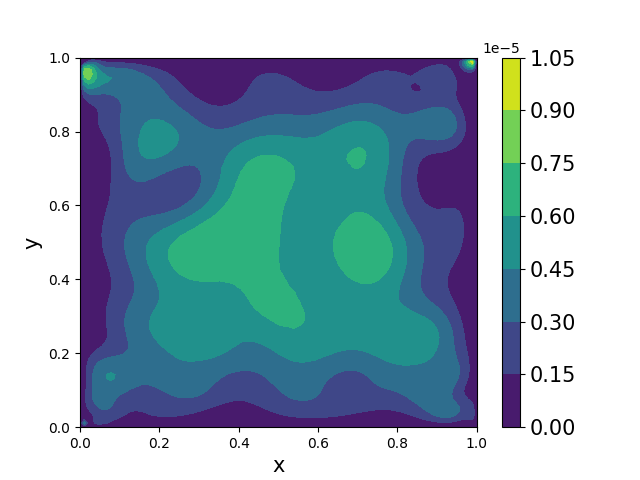}
     \end{subfigure}
        \caption{PINN solution of the 2D Poisson's boundary value problem on the left and the absolute error compared to a converged finite difference solution seen as the ground truth on the right.}
        \label{fig:poisson_sol}
\end{figure}

The PINN solution of the boundary value problem is shown in figure~\ref{fig:poisson_sol} along with the absolute error calculated using the finite difference solution.

The sensitivity terms $\frac{\partial u}{\partial k_i}$ can then be plotted to see the effect of the diffusivity on the solution.

\begin{figure}[H]
    \centering
    \includegraphics[width=0.8\textwidth]{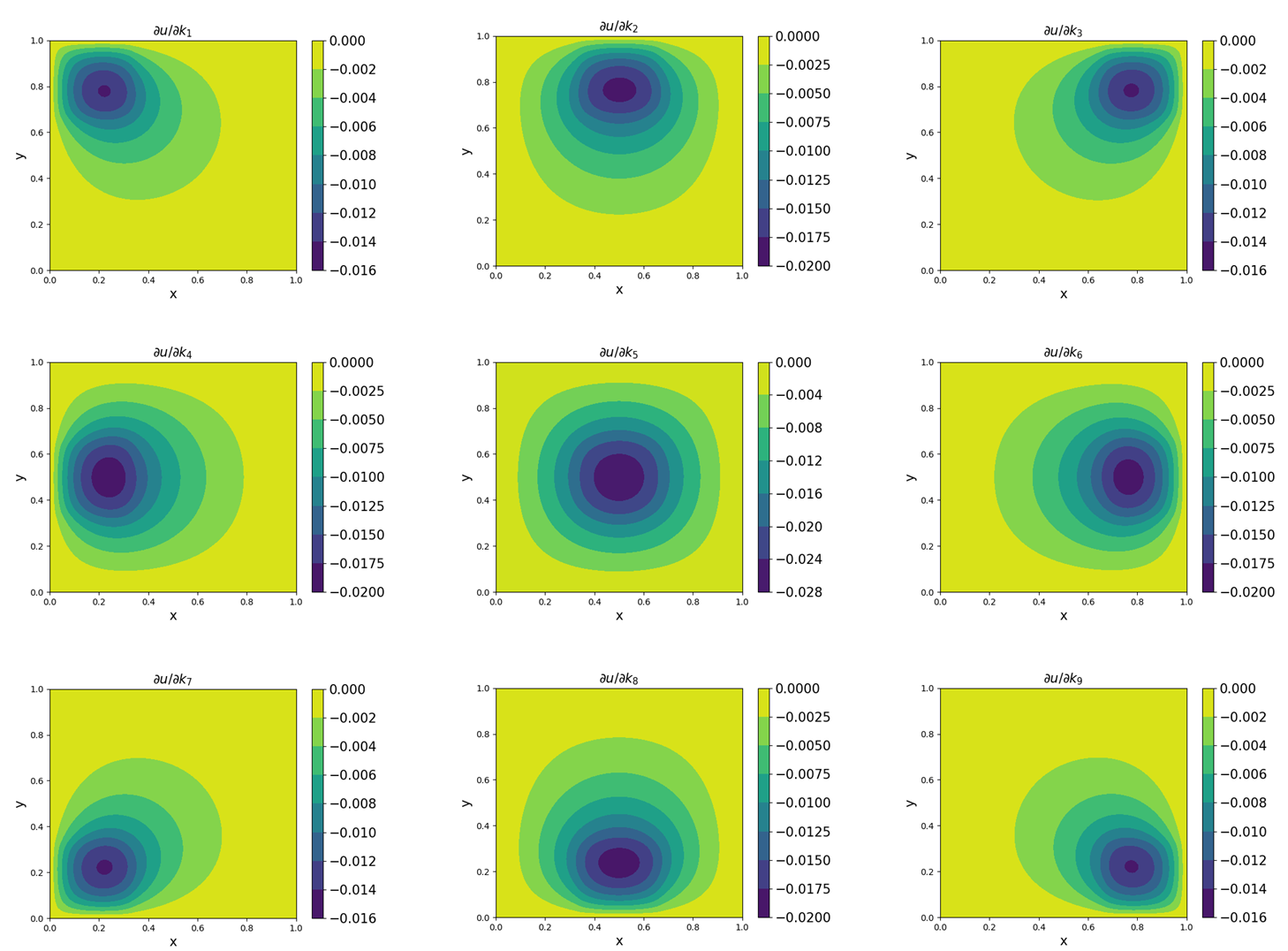}
  \caption{Derivatives of the solution with respect to $k_i$ of the 2D Poisson's problem.}
  \label{fig:poisson_ders}
\end{figure}

\begin{figure}[H]
    \centering
    \includegraphics[width=0.8\textwidth]{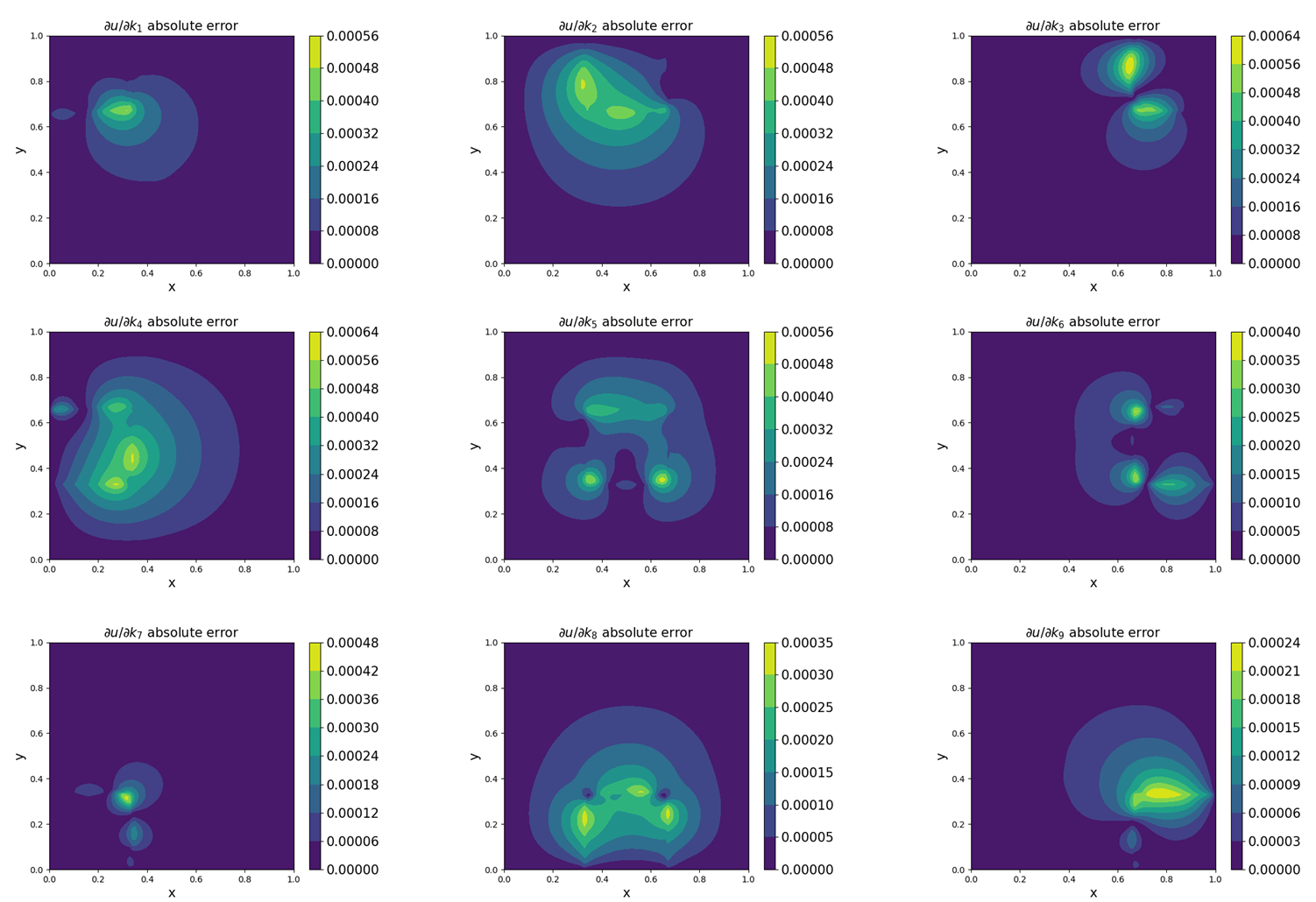}
  \caption{Absolute error of the different derivatives of the solution with respect to $k_i$ of the 2D Poisson's problem.}
  \label{fig:poisson_ders_errors}
\end{figure}

The absolute errors of the calculated derivative are obtained using the finite difference converged solution. The error plots are shown in figure~\ref{fig:poisson_ders_errors}.

\begin{figure}[H]
    \centering
    \includegraphics[width=0.6\textwidth]{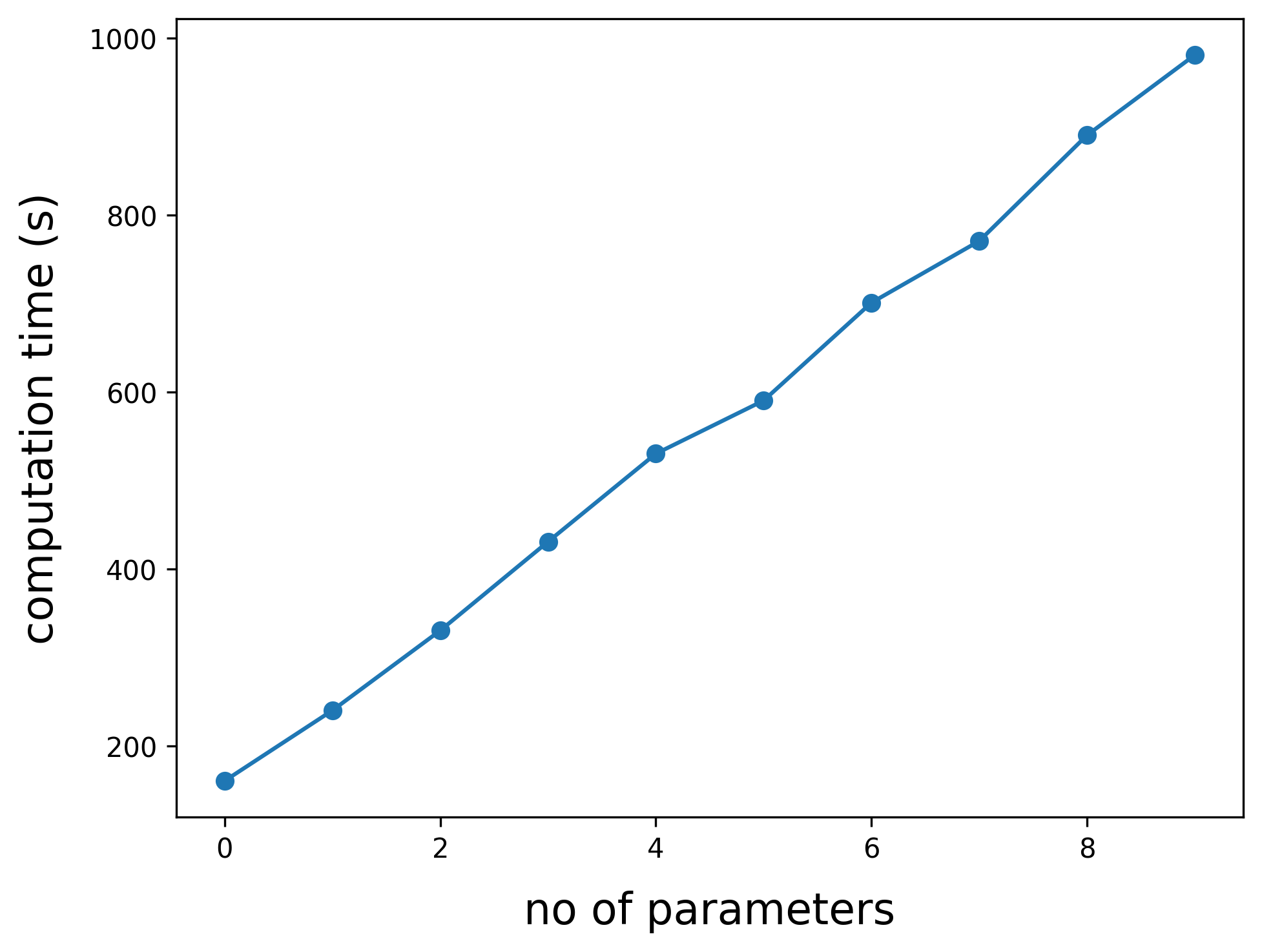}
  \caption{Computational time vs. no of sensitivity parameters.}
  \label{fig:time_poisson}
\end{figure}

The computational time is plotted versus the number of parameters with respect to which sensitivity terms are added in figure~\ref{fig:time_poisson}.

It can be seen from the figure that the computational time grows linearly when increasing the number of parameters the sensitivity is calculated with respect to. This happens because the number of collocation points is the same when adding a new term to the loss function; the added cost is the same when adding new sensitivity terms.  


\subsection{1D two-phase flow in porous media}

\begin{figure}[H]
    \centering

\tikzset{every picture/.style={line width=0.75pt}} 
\begin{tikzpicture}[x=0.75pt,y=0.75pt,yscale=-1,xscale=1]

\draw  [fill={rgb, 255:red, 128; green, 128; blue, 128 }  ,fill opacity=1 ] (123,93) -- (244.8,93) -- (244.8,133) -- (123,133) -- cycle ;
\draw   (244.8,93) -- (535.8,93) -- (535.8,133) -- (244.8,133) -- cycle ;
\draw    (77.8,112.8) -- (119.8,113.75) ;
\draw [shift={(121.8,113.8)}, rotate = 181.3] [color={rgb, 255:red, 0; green, 0; blue, 0 }  ][line width=0.75]    (10.93,-3.29) .. controls (6.95,-1.4) and (3.31,-0.3) .. (0,0) .. controls (3.31,0.3) and (6.95,1.4) .. (10.93,3.29)   ;
\draw    (227.8,167.8) .. controls (267.4,138.1) and (203.11,153.48) .. (240.63,124.69) ;
\draw [shift={(241.8,123.8)}, rotate = 503.13] [color={rgb, 255:red, 0; green, 0; blue, 0 }  ][line width=0.75]    (10.93,-3.29) .. controls (6.95,-1.4) and (3.31,-0.3) .. (0,0) .. controls (3.31,0.3) and (6.95,1.4) .. (10.93,3.29)   ;
\draw    (211.8,66.8) .. controls (216.75,81.65) and (229.44,97.48) .. (185.26,112.54) ;
\draw [shift={(183.9,113)}, rotate = 341.64] [color={rgb, 255:red, 0; green, 0; blue, 0 }  ][line width=0.75]    (10.93,-3.29) .. controls (6.95,-1.4) and (3.31,-0.3) .. (0,0) .. controls (3.31,0.3) and (6.95,1.4) .. (10.93,3.29)   ;
\draw    (404.8,69.8) .. controls (409.75,84.65) and (422.44,100.48) .. (378.26,115.54) ;
\draw [shift={(376.9,116)}, rotate = 341.64] [color={rgb, 255:red, 0; green, 0; blue, 0 }  ][line width=0.75]    (10.93,-3.29) .. controls (6.95,-1.4) and (3.31,-0.3) .. (0,0) .. controls (3.31,0.3) and (6.95,1.4) .. (10.93,3.29)   ;
\draw (78,66.4) node [anchor=north west][inner sep=0.75pt]  [font=\Large]  {$p_{i}{}_{n}$};
\draw (550,60.4) node [anchor=north west][inner sep=0.75pt]  [font=\Large]  {$p_{o}{}_{u}{}_{t}$};
\draw (141,159) node [anchor=north west][inner sep=0.75pt]   [align=left] {{\large Flow front}};
\draw (185,43) node [anchor=north west][inner sep=0.75pt]   [align=left] {fluid 2};
\draw (379,46) node [anchor=north west][inner sep=0.75pt]   [align=left] {fluid 1};
\end{tikzpicture}
    \caption{One-dimensional domain (filling problem)}
    \label{fig:1d_domain}
\end{figure}
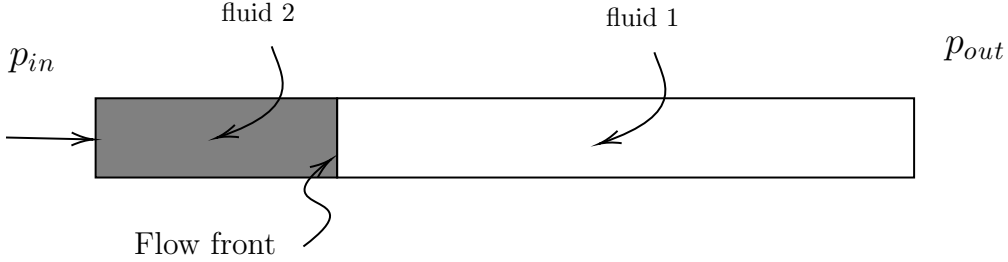

In this section, we introduce a 1D two-phase flow in porous media problem. The problem is faced in liquid composite molding manufacturing processes for instance, where resin is injected into a mold that has a prepositioned fibrous matrix. The problem is shown in figure~\ref{fig:1d_domain}. At $t=0$, the domain is initially saturated with one fluid (fluid 1). Another fluid (fluid 2) is being injected from the left end at constant pressure $p_{in}$, while the pressure at the other end is fixed to $p_{out}$.

The momentum equation can be approximated with Darcy's law that can be written in 1D as follows:

\begin{equation}
    v = -\dfrac{k}{\phi\mu} p_x
    \label{eq:darcy_full}
\end{equation}

where $v$ is the volume average Darcy's velocity, $\mu$ the viscosity, and $p_x$ the pressure gradient, and $\phi$ the porosity. Both fluids are assumed to be incompressible, therefore, the mass conservation equation reduces to

\begin{equation}
   v_x  = 0
    \label{eq:momentum_full}
\end{equation}

Pressure boundary conditions can prescribed on the inlet and outlet:

\begin{equation}\label{eq:p_BC_in}
    p(\mathbf{x_{inlet}},t) = p_{in} ,\ \ p(\mathbf{x_{outlet}},t) = p_{out}
\end{equation}

To track the interface between the two fluids, the Volume Of Fluids (VOF) technique is used; a fraction function $c$ is introduced which takes a value $1$ for the resin and $0$ for the air. The viscosity $\mu$ is redefined as

\begin{equation}
    \mu = c\mu_2 + (1-c)\mu_1
    \label{eq:mu}
\end{equation}

where $\mu_2$ and $\mu_1$ are the two fluids' viscosities. $c$ evolves with time according to the following advection equation

\begin{equation}
    c_t + v c_x  = 0
    \label{eq:c_full}
\end{equation}

where $c_t$ and $c_x$ are the time and spatial derivative of the fraction function $c$, respectively.

Initial and boundary conditions are defined to solve the advection of $c$.

\begin{equation}\label{eq:c_IC}
    c(\mathbf{x},t=0)=c_0(\mathbf{x}),\ \ c(\mathbf{x_{inlet}},t) = 1
\end{equation}

To sum up, the strong form of the problem can be written as:

\begin{equation}
\boxed{
    \begin{gathered}
    c_t + v\ c_x = 0,\ \ \ x \in [0,l],\ \ t\in [0, T],\ \\
    v = -\frac{k}{\phi\mu} p_x,\ \ \ \ \ \ \ x \in [0,l],\ \ t\in [0, T],\ \\
    v_x = 0,\ \ \ \ \ \ x \in [0,l],\ \ t\in [0, T],\ \\
    \mu = c\mu_2 + (1-c)\mu_1\\
    p(0,t) = p_{in},\ \ p(l,t) = p_{out}\\
    c(0,t) = 1,\ \ c(x, 0) = 0
\end{gathered}}
\end{equation}

The parameters of the problem are shown in table~\ref{t1}.

\begin{table}[H]
\centering
\caption{Parameters used for the two-phase flow problem.}
\begin{tabular}{|c|c|}
	\hline
	Parameter & Value \\
	\hline\hline
	$k$ & $1$\\
	\hline
	$\mu_1$ & $10^-5$\\
	\hline
	$\mu_2$ & $1$\\
	\hline
	$p_{in}$ & $1$\\
	\hline
	$p_{out}$ & $0$\\
	\hline
	$l$ & $1$\\
	\hline
    $\phi$ & $1$\\
	\hline
\end{tabular}
\label{t1}
\end{table}

The main PINN terms weights are set to $1$ and $0.1$ for the added sensitivity terms. The adaptivity algorithm presented in \cite{mine} is used to get a better sharper solution reaching 2500 collocation points.

\begin{figure}[H]
    \centering
    \includegraphics[width=0.7\textwidth]{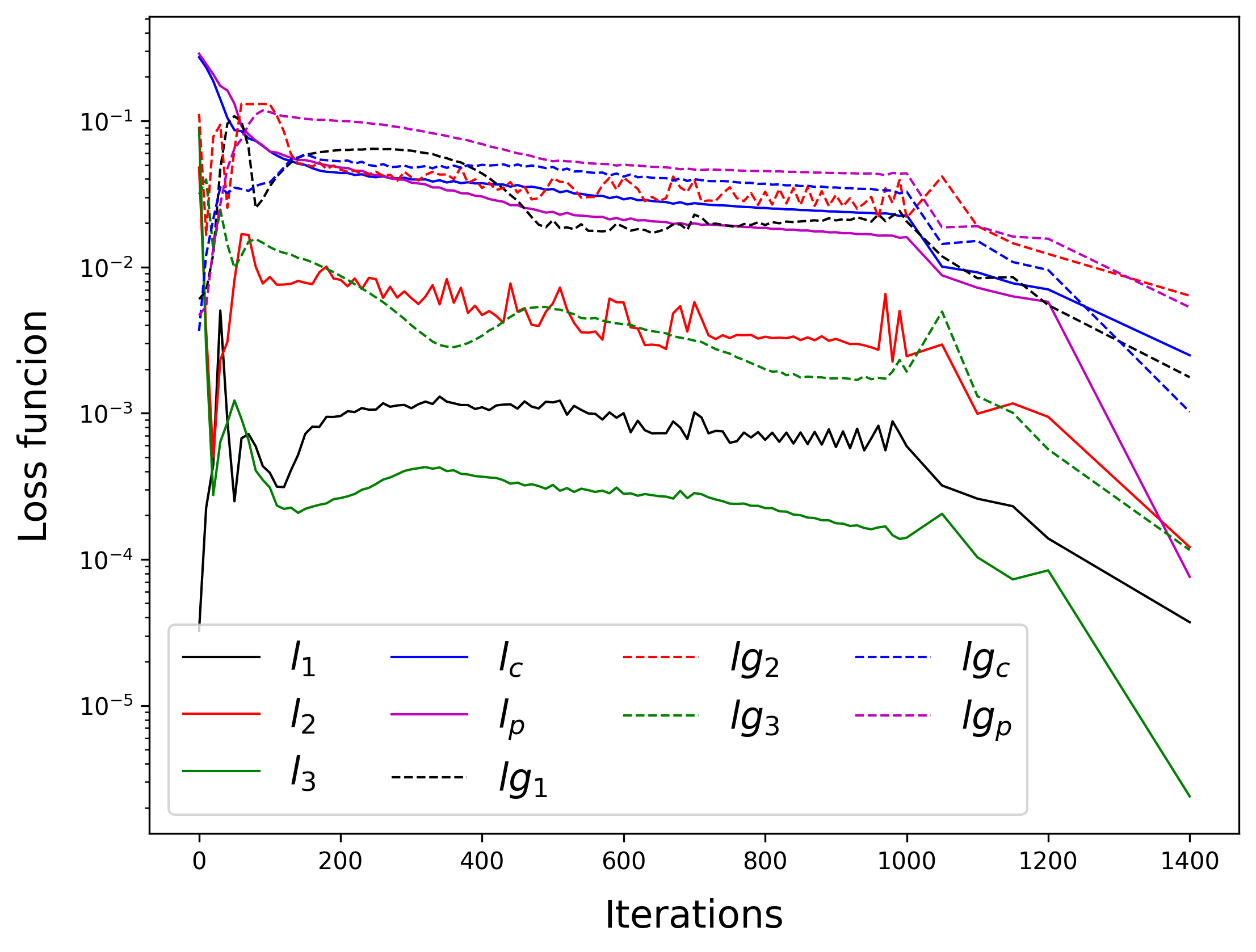}
  \caption{Loss function contributions vs. the number of iterations for the 1D two-phase flow in porous media example.}
  \label{fig:1d_rtm_loss}
\end{figure}

The loss function terms are plotted against the iterations in figure~\ref{fig:1d_rtm_loss}.

\begin{figure}[H]
    \centering
    \includegraphics[width=0.9\textwidth]{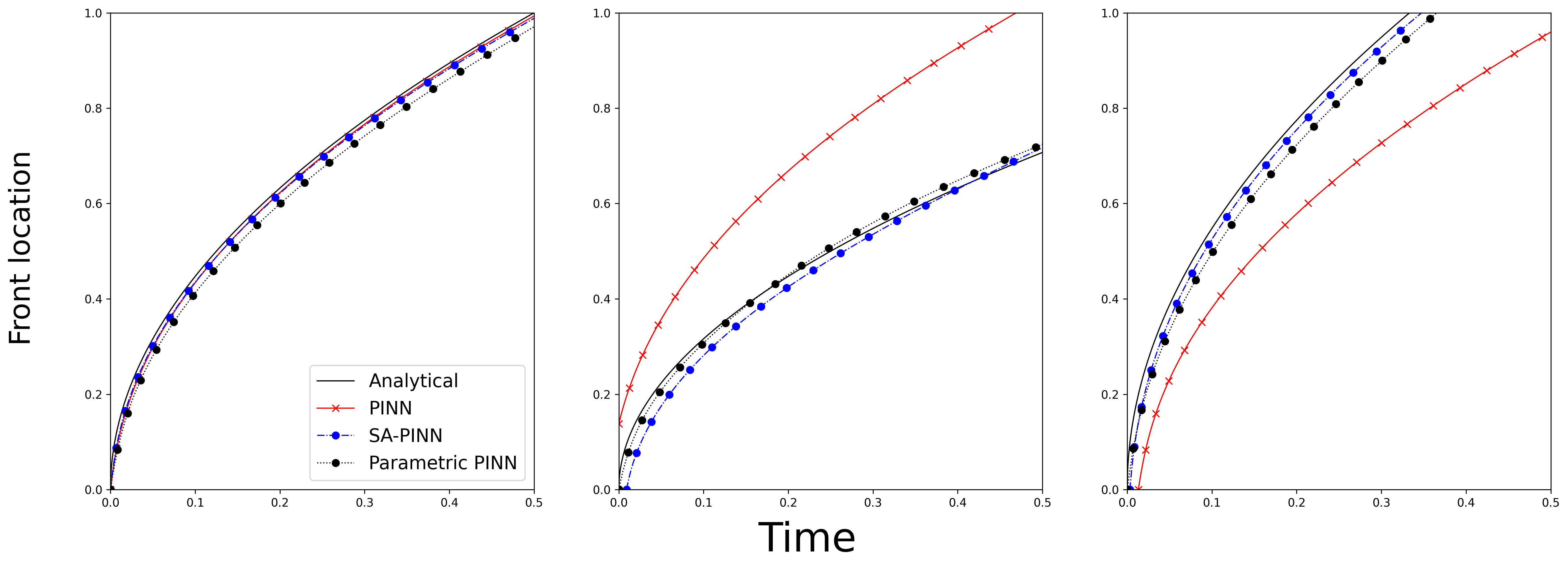}
  \caption{Flow front location vs. time for three different values of $k$ ($k=1$, $0.5$, and $1.5$ from left to right) of the transient two-phase flow in porous media problem.}
  \label{fig:xf}
\end{figure}

First, the front location is plotted for three different values of $k$ by taking the 0.5 level set of the fraction function $c$ in figure~\ref{fig:xf}. A comparison is made between SA-PINN, parametric PINN (collocation points sampled in the parametric space $k\in[0.5, 2]$), and vanilla PINN along with the analytical solution.

It can be noticed that SA-PINN provides good results for values of $k$ away from the nominal value $k=1$. Parametric PINN also succeded in obtaining acceptable results for different values of $k$. However, Vanilla PINN accurately predicts the solution only at the nominal values, however, away from that values, random solutions were obtained which is clear from the two red lines.

For parametric PINN, collocation points were sampled in the parametric space which led to a higher computational cost than both Vanilla and SA-PINN. Using the same machine, the training of vanilla PINN took almost 4 minutes, while SA-PINN nearly 7 minutes. However, parametric PINN took almost 15 minutes to train. This proves that SA-PINN can produce similar results as parametric PINN with much lower computational cost since there are no points sampled in the parametric space in this case.

\begin{figure}[H]
    \centering
    \includegraphics[width=0.7\textwidth]{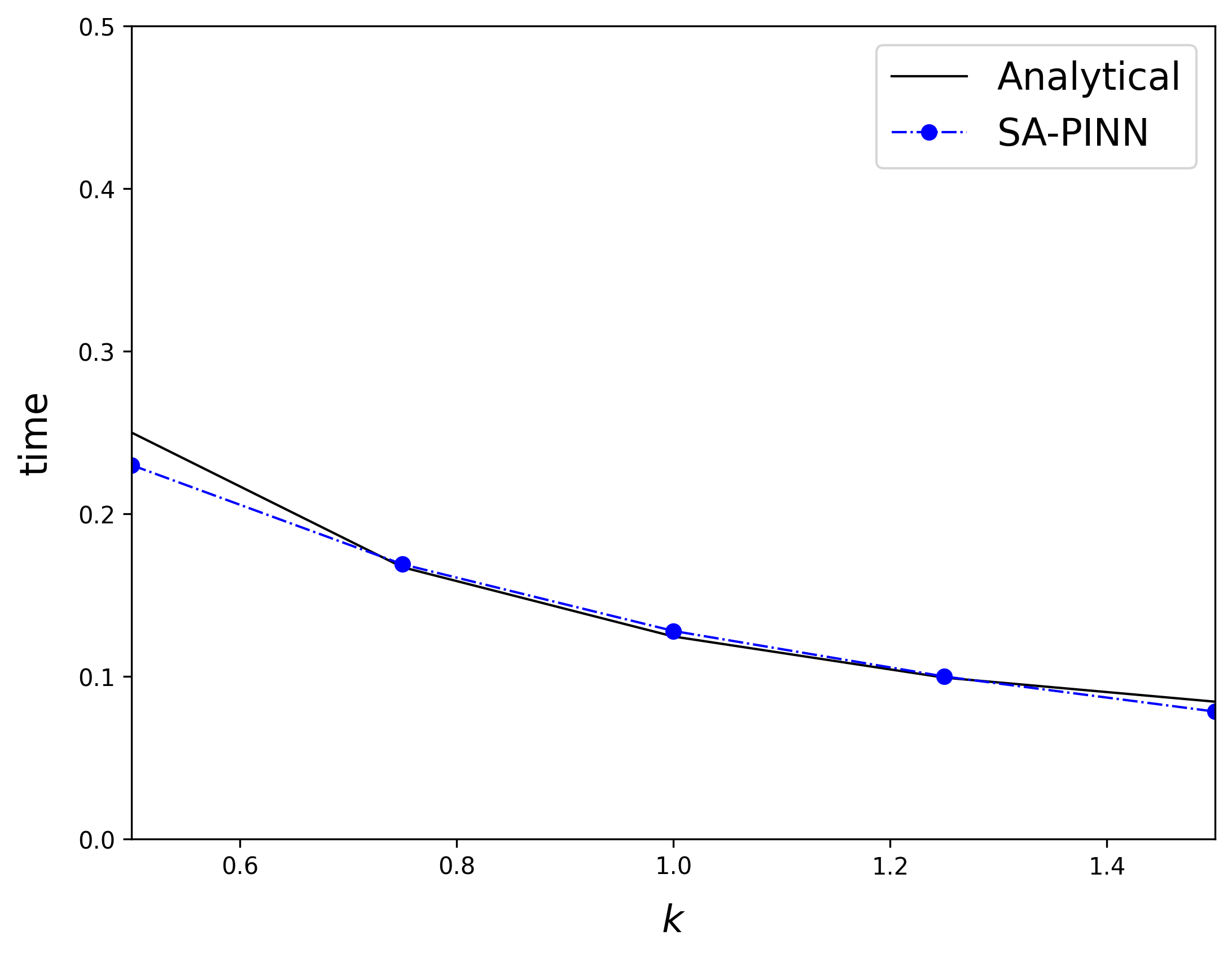}
  \caption{Time at which the flow front reaches $x=0.5$ vs. $k$ for the transient two-phase flow in porous media problem.}
  \label{fig:fill_time}
\end{figure}

In figure~\ref{fig:fill_time}, the time the flow front reaches $x=0.5$ vs. $k$ is plotted. The SA-PINN solution is compared with the analytical solution.

We can see a good estimation of the filling time at different values of $k$ using SA-PINN. This result can be useful in applications of injection processes to estimate the filling time as a function of a parameter of interest.

\subsection{2D two-phase flow in porous media}

\begin{figure}[H]
    \centering

\tikzset{every picture/.style={line width=0.75pt}} 

\begin{tikzpicture}[x=0.75pt,y=0.75pt,yscale=-1,xscale=1]

\draw  [fill={rgb, 255:red, 8; green, 84; blue, 175 }  ,fill opacity=1 ] (227.82,24.32) -- (450,24.32) -- (450,246.5) -- (227.82,246.5) -- cycle ;
\draw    (205.5,222) -- (227.5,222) ;
\draw    (206,246.86) -- (218,246.86) -- (228,246.86) ;
\draw    (178,234.5) -- (211,233.89) ;
\draw [shift={(213,233.86)}, rotate = 178.95] [color={rgb, 255:red, 0; green, 0; blue, 0 }  ][line width=0.75]    (10.93,-3.29) .. controls (6.95,-1.4) and (3.31,-0.3) .. (0,0) .. controls (3.31,0.3) and (6.95,1.4) .. (10.93,3.29)   ;
\draw    (449.76,114) -- (471.76,114) ;
\draw    (449.33,133.86) -- (461.33,133.86) -- (471.33,133.86) ;
\draw    (449.33,123.86) -- (488,123.52) ;
\draw [shift={(490,123.5)}, rotate = 179.5] [color={rgb, 255:red, 0; green, 0; blue, 0 }  ][line width=0.75]    (10.93,-3.29) .. controls (6.95,-1.4) and (3.31,-0.3) .. (0,0) .. controls (3.31,0.3) and (6.95,1.4) .. (10.93,3.29)   ;
\draw    (227.82,234.64) -- (450,233.29) ;
\draw   (227.82,241.39) .. controls (231.89,244.01) and (235.79,246.5) .. (240.32,246.5) .. controls (244.84,246.5) and (248.74,244.01) .. (252.82,241.39) .. controls (256.89,238.78) and (260.79,236.29) .. (265.32,236.29) .. controls (269.84,236.29) and (273.74,238.78) .. (277.82,241.39) .. controls (281.89,244.01) and (285.79,246.5) .. (290.32,246.5) .. controls (294.84,246.5) and (298.74,244.01) .. (302.82,241.39) .. controls (306.89,238.78) and (310.79,236.29) .. (315.32,236.29) .. controls (319.84,236.29) and (323.74,238.78) .. (327.82,241.39) .. controls (331.89,244.01) and (335.79,246.5) .. (340.32,246.5) .. controls (344.84,246.5) and (348.74,244.01) .. (352.82,241.39) .. controls (356.89,238.78) and (360.79,236.29) .. (365.32,236.29) .. controls (369.84,236.29) and (373.74,238.78) .. (377.82,241.39) .. controls (381.89,244.01) and (385.79,246.5) .. (390.32,246.5) .. controls (394.84,246.5) and (398.74,244.01) .. (402.82,241.39) .. controls (406.89,238.78) and (410.79,236.29) .. (415.32,236.29) .. controls (419.84,236.29) and (423.74,238.78) .. (427.82,241.39) .. controls (431.89,244.01) and (435.79,246.5) .. (440.32,246.5) .. controls (443.7,246.5) and (446.74,245.1) .. (449.76,243.31) ;
\draw    (166.82,109.64) .. controls (206.42,79.94) and (248.96,111.98) .. (296.38,131.06) ;
\draw [shift={(297.82,131.64)}, rotate = 201.6] [color={rgb, 255:red, 0; green, 0; blue, 0 }  ][line width=0.75]    (10.93,-3.29) .. controls (6.95,-1.4) and (3.31,-0.3) .. (0,0) .. controls (3.31,0.3) and (6.95,1.4) .. (10.93,3.29)   ;
\draw    (504.82,221.64) .. controls (453.14,197.26) and (426.18,216.62) .. (413.75,236.14) ;
\draw [shift={(412.82,237.64)}, rotate = 300.96] [color={rgb, 255:red, 0; green, 0; blue, 0 }  ][line width=0.75]    (10.93,-3.29) .. controls (6.95,-1.4) and (3.31,-0.3) .. (0,0) .. controls (3.31,0.3) and (6.95,1.4) .. (10.93,3.29)   ;

\draw (149,218) node [anchor=north west][inner sep=0.75pt]    {$v_{in}$};
\draw (496,111) node [anchor=north west][inner sep=0.75pt]    {$p_{out}$};
\draw (147,107) node [anchor=north west][inner sep=0.75pt]    {$k_{1}$};
\draw (508,213) node [anchor=north west][inner sep=0.75pt]    {$k_{2}$};
\end{tikzpicture}
    \caption{Two-dimensional domain (filling problem)}
    \label{fig:2d_filling}
\end{figure}
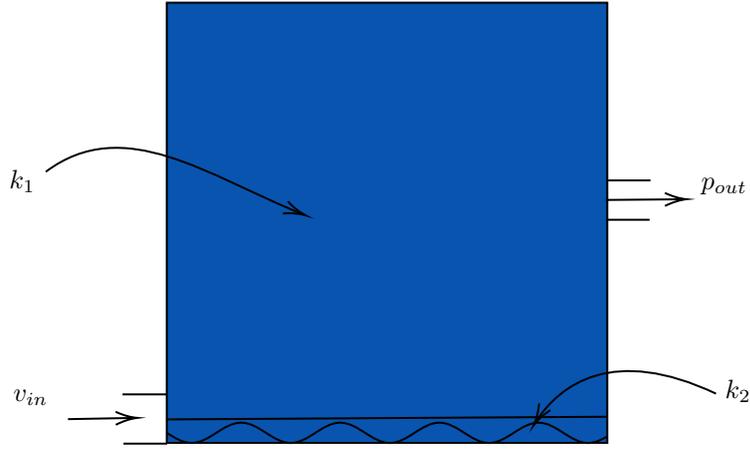

In this example, a 2D two-phase flow in porous media problem is solved which has significance in the liquid composite molding manufacturing process. In the problem, a square mold is under consideration where the inlet exists at the bottom of the left wall and the outlet at the right wall. While the rest of the walls are impermeable. A velocity inlet condition is applied. The whole domain has the same permeability except for a rectangular region at the bottom of the mold where the permeability is significantly higher which represents a common defect in the manufacturing process namely race tracking \cite{race_tracking}. The domain and problem configuration is shown in figure~\ref{fig:2d_filling}.

\begin{figure}[H]
    \centering
    \includegraphics[width=0.7\textwidth]{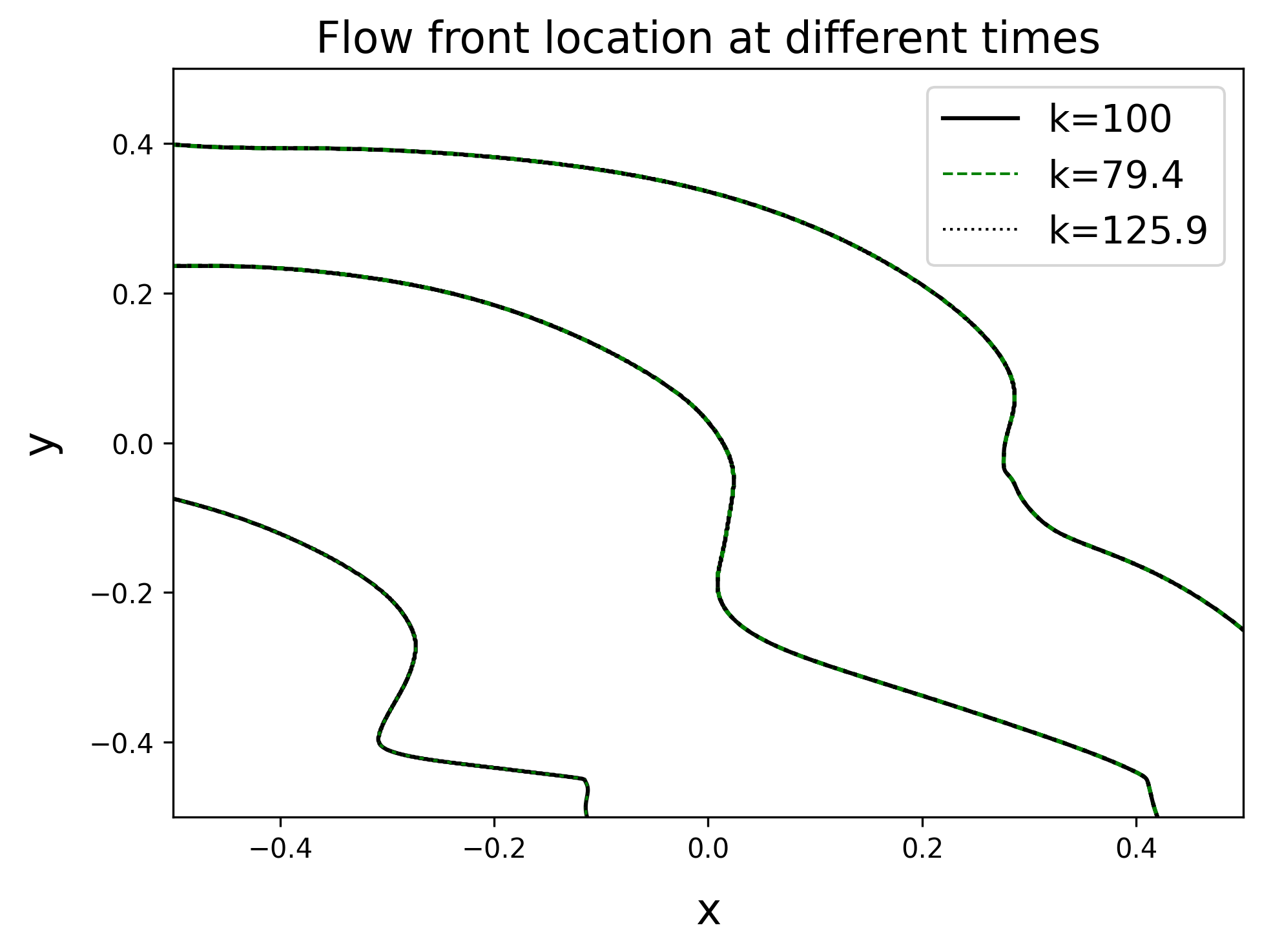}
  \caption{Front locations at different time instants and for different permeability ratios for the 2D two-phase flow in porous media.}
  \label{fig:front_locations}
\end{figure}

The permeability of the whole domain is fixed $k_1=1$. The ratio of permeabilities $k=\frac{k_2}{k_1}$ is chosen to be the parameter of interest for this study where the nominal value is chosen to be 100. At first, the flow front location is plotted at different times and for different permeability ratios in figure~\ref{fig:front_locations}.

It can be seen from figure~\ref{fig:front_locations} that the flow is faster in the higher permeability (race tracking) region which is shown by the advancement of the flow in this region as compared to the bulk of the domain.

\begin{figure}[H]
    \centering
    \includegraphics[width=0.7\textwidth]{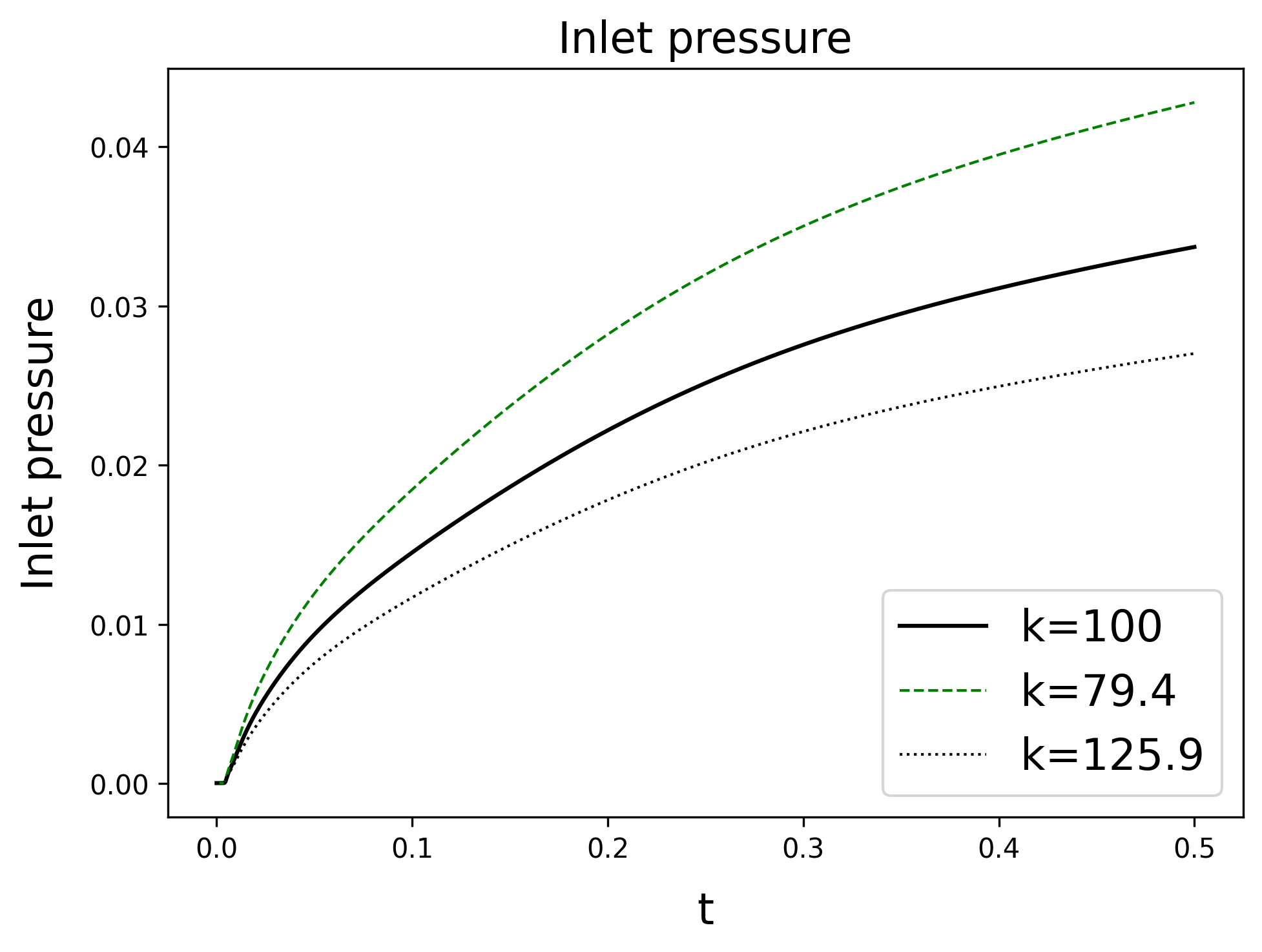}
  \caption{Inlet pressure for different permeability ratios for the 2D two-phase flow in porous media.}
  \label{fig:pressures}
\end{figure}

It can also be noted that the front location evolution is almost the same for different permeability ratios which can be due to the velocity inlet condition that leads to the same flow behavior. Moreover,  the permeability difference is not big enough to cause a noticeable flow behavior change. However, the inlet pressure clearly shows the effect of the permeability ratio change which is shown in figure~\ref{fig:pressures}.

It can be seen from figure~\ref{fig:pressures} that the higher the permeability ratio (high $k_2$), the inlet pressure is reduced. This result physically makes sense as when the permeability gets higher and under constant velocity inlet condition (constant flow rate), less power is needed to drive the flow which leads to less pressure inlet to have the same flow inside the domain. This example has significance in composite manufacturing since it can be used to understand how the defect characteristics (race-tracking) affect the required power to drive the flow in the mold.

\section{Discussion and conclusion}\label{section:conclusion}

The papers introduced a technique to perform sensitivity analysis of input parameters to PDEs using the framework of PINN. The method includes adding the parameters of interest as inputs to the neural network along with the spatial and temporal parameters. Terms representing the derivative of the residual and conditions with respect to the parameters are added to the loss function. Through this combined minimization, we ensure the validity of the solution within a neighborhood of the nominal values of the parameters which allows accurate sensitivity estimation. The method can also be seen as a way to build parametric models, however, it should be treated with care since it can work when the parameter of interest varies only within a small range near a nominal value.

\begin{figure}[H]
    \centering
    \includegraphics[width=0.9\textwidth]{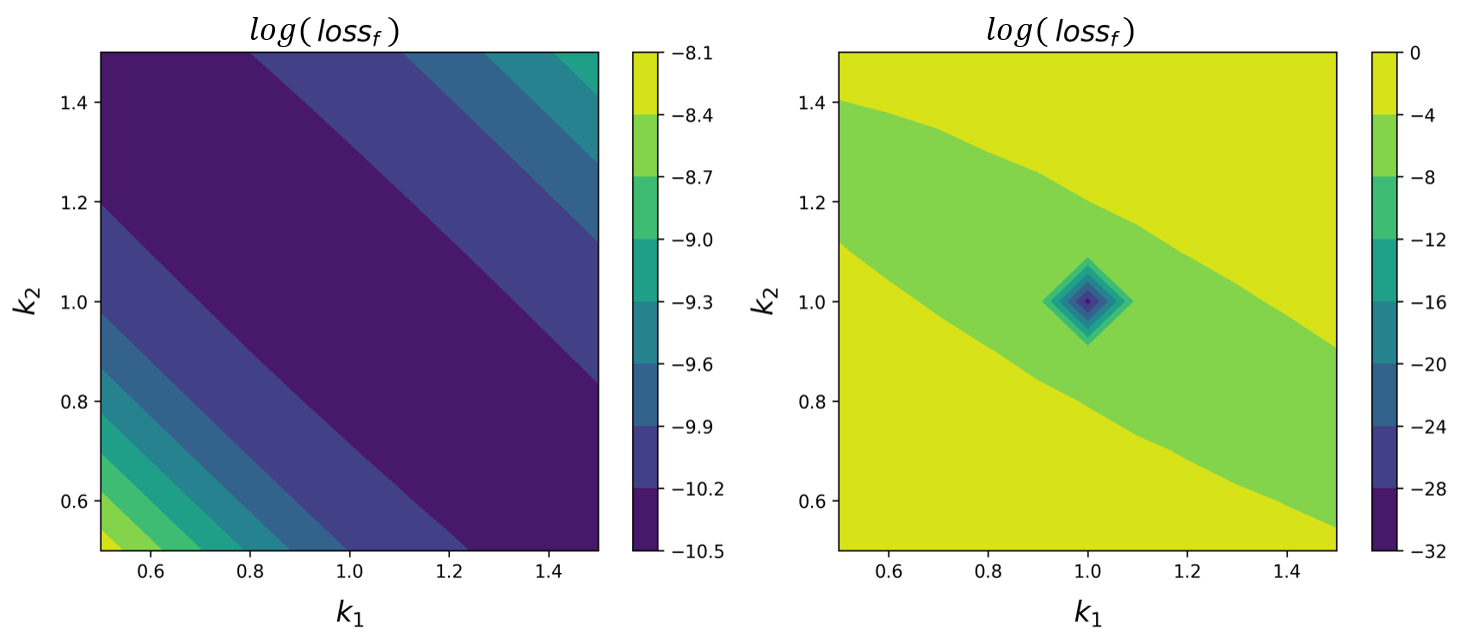}
  \caption{The log contours of $loss_f$ for different $k_1$ and $k_2$ values are plotted of the 2D Poisson's example. On the left, the SA-PINN and on the right the PINN solution. }
  \label{fig:contours}
\end{figure}

In figure~\ref{fig:contours}, we show the effect of the method by plotting the contours of the term $loss_f$ in equation~\ref{eq:loss_2d_poisson} for different $k_i$ values. The plot is done for only 2 parameters for illustration purposes. 

We can see from figure~\ref{fig:contours} that SA-PINN technique has the effect of reducing the main loss function contribution in a neighborhood near the nominal values of the parameters of interest, thus obtaining good solutions in this neighborhood. On the contrary, PINN has an accurate solution only at the nominal values of the parameters. We showed also that the technique does not face problems in dealing with discontinuities as shown in example 3 which could be useful for several applications involving discontinuities.

\section*{Acknowledgements}

This study was funded under the PERFORM Thesis program of IRT Jules Verne, Nantes, France.

\bibliography{biblio}

\end{document}